
  %

  \font\tensc=cmcsc10 
  \font\twelvebf=cmbx10 scaled \magstep1 

  \magnification=\magstep1
  \nopagenumbers
  \voffset=2\baselineskip
  \advance\vsize by -\voffset
  \headline{\ifnum\pageno=1 \hfil\else\ifodd\pageno
    \tensc\hfil fredholm operators and morita equivalence\hfil\folio
    \else
    \tensc\folio\hfil ruy exel\hfil\fi\fi}

  \def\less{<}
  \def\<{\langle}
  \def\>{\rangle}
  \def\:{\colon}
  \def\+{\oplus}
  \def\*{\otimes}
  \def\arw{\rightarrow}
  \def\cstar{$C^*$}
  \def\C{{\bf C}}
  \def\R{{\bf R}}
  \def\N{{\bf N}}
  \def\KK{\hbox{$K\!K$}}
  \def\Im{{\rm Im}}
  \def\Ker{{\rm Ker}}
  \def\sameauthor{\underbar{\hbox to 1.5truecm{\hfil}}}
  \def\square{\hbox{$\sqcap\!\!\!\!\sqcup$}}


\def\statement#1#2{\medbreak\noindent{\bf#1.\enspace}{\sl#2}\medbreak}
  \def\proofbegin{\medbreak\noindent{\bf Proof.\enspace}}
  \def\proofend{\hfill\square\medbreak}
  \def\today{\ifcase\month\or January\or February\or March\or April\or
    May\or June\or July\or August\or September\or October\or
    November\or December\fi\ \number\day, \number\year}
  \def\cite#1{{\rm[\bf #1\rm]}}
  \def\stdbib#1#2#3#4#5#6#7#8{\smallskip \item{[#1]} #2, ``#3'',
    {\sl#4} {\bf#5} (#6), #7--#8.}
  \def\bib#1#2#3#4{\smallskip \item{[#1]} #2, ``#3'', {#4.}}

  \def\larw{\longrightarrow}
  \def\At{\tilde A}
  \def\K{{\cal K}}
  \def\B{{\cal B}}
  \def\FA{{\cal F}_A}
  \def\KA{{\cal K}_A}
  \def\LA{{\cal L}_A}
  \def\LAt{{\cal L}_{\At}}
  \def\LB{{\cal L}_B}
  
  \def\Tt{\tilde{T}}
  \def\St{\tilde{S}}
  \def\ep{\varepsilon}
  \def\mat#1#2#3#4{\pmatrix{#1 & #2 \cr #3 & #4}}
  \def\rank{{\rm rank}}
  \def\ind{{\rm ind}}

  \def\Blackadar {1}
  \def\Brown {2}
  \def\BGR {3}
  \def\BMS {4}
  \def\Connes {5}
  \def\Jensen {6}
  \def\KI {7}
  \def\KII {8}
  \def\Paschke {9}
  \def\RI {10}
  \def\RII {11}
  \def\RIII {12}
  \def\Rowen {13}

  \null
  \vskip -\voffset
  \vskip -1.5truecm
  \baselineskip=9pt
  \rightline{\sevenrm UNM--RE--004}
  \rightline{\sevenrm Concluded December 29, 1992}
  \rightline{\sevenrm Printed \today}

  \vskip 2.5truecm
  \baselineskip=20pt
  \centerline{\twelvebf A FREDHOLM OPERATOR APPROACH TO}
  \centerline{\twelvebf MORITA EQUIVALENCE}

  \bigskip
  \bigskip
  \centerline{\tensc Ruy Exel
  \footnote{\raise.5ex\hbox{*}}
  {\rm On leave from the University of S\~ao Paulo.}}

  \bigskip
  \baselineskip=12pt
  \centerline{\it Department of Mathematics and Statistics}
  \centerline{\it University of New Mexico}
  \centerline{\it Albuquerque, New Mexico 87131}
  \centerline{\it e-mail: exel@math.unm.edu}

  \vfill
  \midinsert\narrower\narrower
  \noindent {\bf Abstract.} Given \cstar-algebras $A$ and $B$ and an
imprimitivity $A$-$B$-bimodule $X$, we construct an explicit
isomorphism $X_*\:K_i(A) \arw K_i(B)$ where $K_i$ denote the complex
$K$-theory functors for $i=0, 1$.  Our techniques do not require
separability nor existence of countable approximate identities. We
thus extend, to general \cstar-algebras, the result of Brown, Green
and Rieffel according to which strongly Morita equivalent
\cstar-algebras have isomorphic $K$-groups. The method employed
includes a study of Fredholm operators on Hilbert modules.
  \endinsert\vfill\eject

\beginsection 1. Introduction

Strong Morita equivalence for \cstar-algebras was introduced by
Rieffel (\cite{\RI}, \cite{\RII}, \cite{\RIII}), generalizing the
concept of Morita equivalence for rings (see for example,
\cite{\Rowen}, Chapter 4).  In \cite{\BGR}, Brown, Green and Rieffel
proved that, if the strongly Morita equivalent \cstar-algebras $A$ and
$B$ are $\sigma$-unital (ie., possess countable approximate
identities), then they are stably isomorphic in the sense that $\K\*A$
is isomorphic to $\K\*B$ ($\K$ denoting the algebra of compact
operators on a separable Hilbert space).  This result has, among many
important consequences, the Corollary that strongly Morita equivalent
$\sigma$-unital \cstar-algebras have isomorphic $K$-theory groups.

In the absence of the $\sigma$-unital hypothesis there are examples of
strongly Morita equivalent \cstar-algebras which are not stably
isomorphic, even if one allows tensoring by the compact operators on a
non-separable Hilbert space (see \cite{\BGR}, Theorem 2.7).  Thus,
without assuming countable approximate identities, the question of
whether the $K$-groups of strongly Morita equivalent \cstar-algebras
are isomorphic, remained open.

Even though the argument can be made that the majority of
\cstar-algebras of interest are $\sigma$-unital, there is a more
serious obstacle in applying the BGR Theorem to concrete situations.
The proof of that result rests largely, on the work of Brown on full
hereditary subalgebras of \cstar-algebras \cite{\Brown}, where a
highly non-constructive method is used.  The end result is that, given
strongly Morita equivalent, $\sigma$-unital \cstar-algebras $A$ and
$B$, one knows that $\K\*A$ and $\K\*B$ are isomorphic to each other
but there is generally little hope that one can exhibit a concrete
isomorphism, unless, of course, an isomorphism is known to exist
independently of the Morita equivalence.  Likewise, no expression for
the isomorphism of $K$-groups may be described in general.

The primary goal of the present work is to make explicit the
isomorphism between $K_*(A)$ and $K_*(B)$ that is predicted by the BGR
Theorem. Recently, Rieffel suggested that, if such an explicit
construction was possible, then it would probably work without
separability, and in fact it does: in Theorems 5.3 and 5.5 we prove
that strongly Morita equivalent \cstar-algebras have isomorphic
$K$-groups, irrespective of countable approximate identities.

It should be noted that in a crucial point of our argument (see the
proof of Theorem 2.7 below), we have to use Kasparov's stabilization
Theorem \cite{\KI} which requires a certain countability condition,
but in Lemma 2.6, we manage to conform ourselves to Kasparov's
hypothesis, without having to assume any countability beforehand.

Let us now briefly describe the strategy adopted to prove our main
result. Given a \cstar-algebra $A$, we begin by studying Fredholm
operators on Hilbert $A$-modules and we define the Fredholm index (3.4
and 3.10) which takes values in $K_0(A)$.  We then prove that any
element of $K_0(A)$ is the index of some Fredholm operator (3.14),
showing the surjectivity of the Fredholm index.  Next, we characterize
in a somewhat geometric way, the pairs of Fredholm operators having
identical index (3.16).  This enables us to introduce an equivalence
relation on the set of all Fredholm operators, under which two
operators are equivalent if and only if they have the same index.  The
quotient by this relation therefore provides an alternate definition
for $K_0(A)$ (3.17).

If one is given \cstar-algebras $A$ and $B$ and a Hilbert
$A$-$B$-bimodule $X$ then, for every Fredholm operator $T\:M\arw N$,
where $M$ and $N$ are Hilbert $A$-modules, one can consider the
operator
  $T\*I_X\: M\*_BX \arw N\*_BX$.
  The crux of the whole matter is to prove that this operator is a
Fredholm operator (4.3) if $X$ is left-full (in the sense that the
range of the $A$-valued inner-product generates $A$, see 4.1).  This
tensor product construction is then used to define (5.1) a group
homomorphism $X_*\:K_0(A)\arw K_0(B)$ satisfying $X_*(\ind(T)) =
\ind(T\*I_X)$.  When $X$ is also right-full we show that $X_*$ is an
isomorphism (5.3). The case of $K_1$ is treated, as usual, using
suspensions (5.5).

Our characterization of $K_0(A)$ in terms of Fredholm operators should
be compared to the well known fact that $K_0(A)$ is isomorphic to the
Kasparov \cite{\KII} group $\KK(\C, A)$ (\C\ denoting the algebra of
complex numbers). However, since we do not assume countable
approximate identities, that fact does not seem to follow from the
existing machinery of \KK-theory.  Furthermore, this aspect of our
work is, perhaps, an indication that \KK-theory may be extended beyond
the realm of $\sigma$-unital \cstar-algebras. In particular, it sounds
reasonable to conjecture that strongly Morita equivalent
\cstar-algebras are \KK-equivalent. In fact, our methods seem
particularly well suited to prove such a conjecture if only \KK-theory
were not so inextricably linked to $\sigma$-unital \cstar-algebras.

To some extent, one should also compare what we do here with the
Fredholm modules of Connes \cite{\Connes}, in the sense that we
explore, in more detail than usually found in the existing literature,
what happens when one plugs the algebra of complex numbers as one of
the variables of Kasparov's \KK-functor. Connes' Fredholm modules are
literally the ingredients of $\KK(\cdot\thinspace,\C)$, according to
the Fredholm picture of \KK-theory (see \cite{\Blackadar}). On the
other hand, the Fredholm operators we study here are related to
$\KK(\C,\cdot)$.

Although Hilbert modules permeate all of our work, we believe the
subject is old enough that we do not need to spend much time
presenting formal definitions from scratch. The reader will find all
of the relevant definitions in \cite{\Paschke}, \cite{\RI},
\cite{\Blackadar}, \cite{\Jensen} and \cite{\BMS}. Nevertheless let us
stress that all of our Hilbert modules are supposed to be right
modules except, of course, when they are bimodules. The term {\it
operator}, when referring to a map $T\:M\arw N$ between the Hilbert
$A$-modules $M$ and $N$, will always mean an element of $\LA(M, N)$,
that is, $T$ should be an adjointable linear map in the sense that
there exists $T^*\:N\arw M$ satisfying $\<T(\mu),\nu\> = \<\mu,
T^*(\nu)\>$ for $\mu\in M$ and $\nu\in N$. See \cite{\Jensen}, 1.1.7
for more details.

\beginsection 2. Preliminaries on Hilbert Modules

In this section we would like to present a few simple facts about
Hilbert modules which we shall need in the sequel.  Throughout this
section, $A$ will denote a fixed \cstar-algebra and $M$ and $N$ will
always refer to (right) Hilbert modules over $A$. As usual, $A^n$ will
be viewed as a Hilbert module equipped with the inner-product
  $\<(a_i)_i, (b_i)_i\> = \sum_{i=1}^n a_i^* b_i$.

Let us observe that $A^n$ will often stand for the set of $n\times 1$
(column) matrices over $A$. In that way, the above inner-product can
be expressed, for $v=(a_i)_i$ and $w=(b_i)_i$, as $\<v, w\> = v^*w$.
Note that $v^*$ refers to the conjugate-transpose matrix.

For each n-tuple $\mu = (\mu_i)_i$ in $M^n$, we denote by $\Omega_\mu$
the operator in $\LA(A^n, M)$ defined by
  $$\Omega_\mu\bigl((a_i)_i\bigr) = \sum_{i=1}^n \mu_i a_i, \quad
(a_i)_i \in A^n.$$
  It is easy to see that $\Omega_\mu^*$ is given by
  $$\Omega_\mu^*(\xi) = \bigl(\<\mu_i,\xi\>\bigr)_i, \quad \xi \in
M.$$ If $\nu = (\nu_i)_i$ is an n-tuple of elements of $N$, then the
operator $T = \Omega_\nu \Omega_\mu^*$ is in $\LA(M, N)$. More
explicitly we have
  $$T(\xi) = \sum_{i=1}^n \nu_i\<\mu_i,\xi\>, \quad \xi \in M.$$
  Maps such as $T$ will be called {\it $A$-finite rank\/} operators
and the set of all those will be denoted $\FA(M, N)$, or just $\FA(M)$
in case $M = N$. The closure of $\FA(M, N)$ in $\LA(M, N)$ is denoted
$\KA(M, N)$ and elements from this set will be referred to as {\it
$A$-compact\/} operators. An expository treatment of operators on
Hilbert modules may be found in \cite{\Jensen}.

\statement{2.1.~Proposition}{For each $\mu=(\mu_i)_i$ in $M^n$ one has
that $\Omega_\mu$ is in $\KA(A^n, M)$ and hence also that
$\Omega_\mu^*$ is in $\KA(M, A^n)$.}

\proofbegin It is obviously enough to consider the case $n=1$.  Let
$(u_\lambda)_\lambda$ be an approximate identity for $A$ (always
assumed to be positive and of norm one). It follows from
\cite{\Jensen}, 1.1.4 that
  $\lim_\lambda \mu u_\lambda = \mu$. Therefore we have for all $a$ in
$A$
  $$\Omega_\mu(a) = \mu a = \lim_\lambda \mu u_\lambda a
  = \lim_\lambda \mu \<u_\lambda, a\>
  = \lim_\lambda \Omega_\mu\Omega_{u_\lambda}^*(a), $$
  the limit being uniform in $\|a\| \leq 1$. \proofend

\statement{2.2.~Definition}{A Hilbert module $M$ will be said to be an
{\it $A$-finite rank} module if the identity operator $I_M$ is in
$\KA(M)$.}

Since $\FA(M)$ is an ideal in $\LA(M)$, which is dense in $\KA(M)$, it
is easy to see that $I_M$ must, in fact, be in $\FA(M)$ whenever $M$
is $A$-finite rank.  We next give the complete characterization of
$A$-finite rank modules.

\statement{2.3.~Proposition}{$M$ is $A$-finite rank if and only if
there exists an idempotent matrix $p$ in $M_n(A)$ such that $M$ is
isomorphic, as Hilbert modules, to $pA^n$.}

\proofbegin Initially we should observe that the use of the term
``isomorphic'', when referring to Hilbert modules, is in accordance
with \cite{\Jensen}, 1.1.18. That is, there should exist a linear
bijection, preserving the $A$-valued inner product.

Assume $M$ to be $A$-finite rank. Then $I_M = \Omega_\nu \Omega_\mu^*$
where $\mu$ and $\nu$ are in $M^n$. Observe that $\Omega_\mu^*
\Omega_\nu$ is then an idempotent $A$-module operator on $A^n$, which
therefore corresponds to left multiplication by the idempotent
$n\times n$ matrix $p = \bigl(\<\mu_i,\nu_j\>\bigr)_{i, j}$. The
operator $\Omega_\mu^*$ then gives an invertible operator in $\LA(M,
pA^n)$.  To make that map an (isometric) isomorphism one uses polar
decomposition. The converse statement is trivial.  \proofend

Any $A$-finite rank module $M$ clearly becomes a finitely generated
projective module over the unitized \cstar-algebra $\At$ (the unitized
algebra is given a new identity element, even if $A$ already has one).

\statement{2.4.~Definition}{The $K$-theory class $[M]_0\in K_0(\At)$
of an $A$-finite rank module $M$, is obviously an element of $K_0(A)$
and will henceforth be denoted $\rank(M)$. If $M$ is not necessarily
assumed to be $A$-finite rank, but if $P$ is an idempotent operator in
$\KA(M)$, then $\Im(P)$ is clearly an $A$-finite rank module. In this
case we let $\rank(P) = \rank(\Im(P))$.}

If $X$, $Y$, $Z$ and $W$ are Hilbert $A$-modules and $T$ is in
$\LA(X\+ Y, Z\+ W)$, then $T$ can be represented by a matrix
  $$T = \mat{T_{ZX}}{T_{ZY}}{T_{WX}}{T_{WY}}, $$
  where $T_{ZX}$ is in $\LA(X, Z)$ and similarly for the other matrix
entries. Matrix notation is used to define our next important concept.

\statement{2.5.~Definition}{The Hilbert modules $M$ and $N$ are said
to be {\it quasi-stably-isomorphic} if there exists a Hilbert module
$X$ and an invertible operator $T$ in $\LA(M\+ X, N\+ X)$ such that
$I_X - T_{XX}$ is $A$-compact.}

Of course the concept just defined is meant to be a generalization of
the well known concept of stable isomorphism for finitely generated
projective $A$-modules, at least when $A$ is unital.  We shall discuss
shortly, the precise sense in which that generalization takes place.
Before that we need a preparatory result.

\statement{2.6.~Lemma}{Assume $M$ and $N$ are $A$-finite rank modules.
If $M$ and $N$ are quasi-stably-isomorphic then the module $X$
referred to in 2.5 can be taken to be countably generated.}

\proofbegin Let
  $T = \mat{T_{NM}}{T_{NX}}{T_{XM}}{T_{XX}}$
  be an invertible operator in $\LA(M\+ X, N\+ X)$ with $I_X-T_{XX}$
$A$-compact, as in 2.5, and let
  $S = \mat{S_{MN}}{S_{MX}}{S_{XN}}{S_{XX}}$
  be the inverse of $T$. Choose a countable set $\Xi_0 =
\{\xi_i\}_{i\in \N}$ of elements in $X$ such that

  \medskip \item{(i)} The images of $T_{XM}$, $S_{XN}$, $T_{NX}^*$,
and $S_{MX}^*$ are contained in the submodule of $X$ generated by
$\Xi_0$.
  \medskip \item{(ii)} $I_X - T_{XX}$ can be approximated by
$A$-finite rank operators of the form $\Omega_\nu \Omega_\mu^*$, where
the components of $\mu = (\mu_1, \ldots, \mu_n)$ and $\nu = (\nu_1,
\ldots, \nu_n)$ belong to $\Xi_0$.

\medskip Define, inductively, $\Xi_{n+1} = \Xi_n \cup T_{XX}(\Xi_n)
\cup S_{XX}(\Xi_n) \cup T_{XX}^*(\Xi_n) \cup S_{XX}^*(\Xi_n)$. The set
$\Xi = \bigcup_{n\in \N} \Xi_n$ is then obviously countable, satisfies
(i) and (ii) above and, in addition,
  \medskip
  \item{(iii)} $\Xi$ is invariant under $T_{XX}$, $S_{XX}$, $T_{XX}^*$
and $S_{XX}^*$.

\medskip Let $X_0$ be the Hilbert submodule of $X$ generated by $\Xi$.
Because of (i) and (iii) we see that $T(M\+ X_0) \subseteq N\+ X_0$
and $T^*(N\+ X_0) \subseteq M\+ X_0$. The restriction of $T$ then
gives an operator $T'$ in $\LA(M\+X_0, N\+X_0)$. The same reasoning
applies to $S$ providing $S'$ in $\LA(N\+X_0, M\+X_0)$ which is
obviously the inverse of $T'$. In virtue of (ii) it is clear that $T'$
satisfies the conditions of definition 2.5. \proofend

\statement{2.7.~Theorem}{If $M$ and $N$ are quasi-stably-isomorphic
$A$-finite rank modules, then $\rank(M) = \rank(N)$.}

\proofbegin According to 2.6 let $X$ be a countably generated Hilbert
module and $T$ be an invertible operator in $\LA(M\+ X, N\+ X)$ such
that $I_X - T_{XX}$ is in $\KA(X)$. By Kasparov's stabilization
Theorem \cite{\KI}, 3.2 (see also \cite{\Jensen}, 1.1.24), $X\+H_A$ is
isomorphic to $H_A$, where $H_A$ is the completion of
$\bigoplus_1^\infty A$ (see \cite{\Jensen}, 1.1.6 for a more precise
definition). This said, we may assume, without loss of generality,
that $X=H_A$. Since $M$ is finitely generated as an $A$-module, we
conclude, again by Kasparov's Theorem, that $M\+X$ is isomorphic to
$H_A$. Choose, once and for all, an isomorphism $\varphi \: H_A \arw
M\+X$ and consider the operators $F$ and $G$ on $H_A$ given by the
compositions
  $$F
  \ :\ H_A \buildrel \varphi \over \larw
  M\+X \buildrel T \over \larw
  N\+X \larw
  X \larw
  M\+X \buildrel \varphi^{-1} \over \larw
  H_A$$ and
  $$G
  \ :\ H_A \buildrel \varphi \over \larw
  M\+X \larw
  X \larw
  N\+X \buildrel T^{-1} \over \larw
  M\+X \buildrel \varphi{-1} \over \larw
  H_A, $$
  where the unmarked arrows denote either the canonical inclusion or
the canonical projection. It can be easily seen that both $I_{H_A} -
GF$ and $I_{H_A} - FG$ belong to $\KA(H_A)$. Consider the exact
sequence of \cstar-algebras
  $$0 \larw \KA(H_A) \larw \LA(H_A) \larw \LA(H_A)/\KA(H_A) \larw 0.$$
  Denoting by $\pi$ the quotient map, one sees that $\pi(F)$ and
$\pi(G)$ are each others inverse. Two facts need now be stressed.  The
first one is that the $K$-theory index map
  $$\ind \ :\ K_1\bigl(\LA(H_A)/\KA(H_A)\bigr) \arw
K_0\bigl(\KA(H_A)\bigr)$$
  assigns to the class of $\pi(F)$, the element $\rank(N) - \rank(M)$,
once $\KA(H_A)$ is identified with $\K\* A$ (according to \cite{\KI},
2.4) and $K_0(\K\* A)$ is identified with $K_0(A)$ as usual in
$K$-theory.

The second fact to be pointed out is that, since $T_{XX}$ is an
$A$-compact perturbation of the identity and $M$ and $N$ are
$A$-finite rank modules (and so any operator having either $M$ or $N$
as domain or codomain must be $A$-compact), one concludes that $F$ is
an $A$-compact perturbation of the identity. It follows that $\pi(F) =
1$ and hence that $\pi(F)$ has trivial index. This concludes the
proof.  \proofend

\beginsection 3. Fredholm Operators

We shall now study Fredholm operators between Hilbert modules.  As
before, $A$ will denote a fixed \cstar-algebra and $M$ and $N$, with
or without subscripts, will denote Hilbert A-modules.

\statement{3.1.~Definition}{Let $T$ be in $\LA(M, N)$. Suppose there
is $S$ in $\LA(N, M)$ such that $I_M - ST$ is in $\KA(M)$ and $I_N -
TS$ is in $\KA(N)$. Then $T$ is said to be an {\it $A$-Fredholm}
operator.  In case the algebra $A$ is understood, we shall just say
that $T$ is a {\it Fredholm} operator (but we should keep in mind that
this notion does not coincide with the classical notion of Fredholm
operators).}

As in the classical theory of Fredholm operators, it can be proved
that whenever $T$ is $A$-Fredholm, one can find $S$ in $\LA(N, M)$
such that $I_M - ST$ is in $\FA(M)$ and $I_N - TS$ is in $\FA(N)$.  In
all of our uses of the $A$-Fredholm hypothesis, below, we shall adopt
that characterization.

In the initial part of the present section we shall concentrate on a
special class of operators which we will call regular operators. This
concept is the natural extension, to Hilbert modules, of the notion of
operators on Hilbert spaces having closed image. Contrary to the
Hilbert space case, not all $A$-Fredholm operators have a closed
image.

\statement{3.2.~Definition}{An operator $T$ in $\LA(M, N)$ is said to
be {\it regular} if there is $S$ in $\LA(N, M)$ such that $TST = T$
and $STS = S$. Any operator $S$ having these properties will be called
a {\it pseudo-inverse} of $T$.}

It is easy to see that for any pseudo-inverse $S$ of $T$ one has that
$I_M - ST$ is the projection onto $\Ker(T)$ and that $TS$ is the
projection onto $\Im(T)$. If $T$ is assumed to be regular and
Fredholm, there are, according to the above definitions, operators $S$
and $S'$ such that $I_M - S'T$ and $I_N - TS'$ are $A$-finite rank
and, on the other hand, $TST = T$ and $STS = S$.

Observe that, because
  $I_M - ST = (I_M - S'T)(I_M - ST)$,
  any pseudo inverse $S$, must be such that $I_M - ST$ is $A$-finite
rank and similarly for $I_N - TS$.  An immediate consequence of the
present discussion is the following:

\statement{3.3.~Proposition}{Let $T \in \LA(M, N)$ be a regular
$A$-Fredholm operator. Then both $\Ker(T)$ and $\Ker(T^*)$ are
$A$-finite rank modules.}

Let us now define the Fredholm index for regular $A$-Fredholm
operators.  Shortly we shall extend that concept to general
$A$-Fredholm operators.

\statement{3.4.~Definition}{If $T$ is a regular $A$-Fredholm operator,
then the {\it Fredholm index} of $T$ is defined to be the element of
$K_0(A)$ given by
  $$\ind(T) = \rank(\Ker(T)) - \rank(\Ker(T^*)).$$}

We collect, in our next proposition, some of the elementary properties
of the Fredholm index.

\statement{3.5.~Proposition}{If $T \in \LA(M, N)$ is a regular
Fredholm operator, then
  \medskip \item{(i)} $\ind(T^*) = -\ind(T)$.
  \medskip \item{(ii)} For any pseudo-inverse $S$ of $T$ we have that
$\rank(\Ker(T^*)) = \rank(\Ker(S))$ and $\ind(S) = -\ind(T)$.
  \medskip \item{(iii)} If $X$ and $Y$ are Hilbert $A$-modules and if
$U \in \LA(X, M)$ and $V \in \LA(N, Y)$ are invertible, then
$\ind(VTU) = \ind(T)$.
  \medskip \item{(iv)} If $T_1 \in \LA(M_1, N_1)$ is regular and
Fredholm, then $\ind(T\+T_1) = \ind(T) + \ind(T_1)$.}

\proofbegin Left to the reader. \proofend

The first fact about classical Fredholm operators whose generalization
to Hilbert modules requires some work is the invariance under compact
perturbations which we now prove.

\statement{3.6.~Theorem}{If $T$ is a regular Fredholm operator in
$\LA(M)$ and if $I_M - T$ is $A$-compact, then $\ind(T) = 0$.}

\proofbegin Let $S$ be a pseudo-inverse for $T$ and denote by $X$ the
image of the idempotent $ST$ or, equivalently, $X= \Im(S)$. Consider
the transformation
  $$U \: \Ker(T) \+ X \arw \Ker(S) \+ X$$
  given by $U(\xi,\eta) = \bigl((I_N - TS)(\xi + \eta), S(\xi +
\eta)\bigr)$. It is easy to see that $U$ is invertible with inverse
given by $V(\xi,\eta) = \bigl((I_M - ST)(\xi + T(\eta)), ST(\xi +
T(\eta))\bigr)$. The operator $U_{XX}$ (occurring in the matrix
representation of $U$) coincides with $S$ which is easily seen to be
an $A$-compact perturbation of the identity. This shows that the
$A$-finite rank modules $\Ker(T)$ and $\Ker(S)$ are
quasi-stably-isomorphic. By 2.7 we conclude that $\rank(\Ker(T)) =
\rank(\Ker(S))$ and hence that $\ind(T) = 0$.  \proofend

\statement{3.7.~Corollary}{If $T_1, T_2 \in \LA(M, N)$ are regular
Fredholm operators such that $T_1 - T_2$ is in $\KA(M, N)$, then
$\ind(T_1) = \ind(T_2)$.}

\proofbegin Let $S_1$ and $S_2$ be pseudo inverses for $T_1$ and
$T_2$. Define operators $U$ and $R$ in $\LA(M\+N)$ by
  $$U = \mat{I_M - S_1 T_1}{S_1}{T_1}{I_N - T_1 S_1}
  \quad {\rm and} \quad
  R = \mat{0}{S_1}{T_2}{0}.$$
  Observe that $U^2 = I$ so, in particular, $U$ is invertible.
Therefore, by 3.5, we have that $\ind(UR) = \ind(R) = \ind(T_2) -
\ind(T_1)$.  But, since
  $$UR = \mat{S_1 T_2} {0} {T_2 - T_1 S_1 T_2} {T_1 S_1}$$
  is an $A$-compact perturbation of the identity, it follows that
$\ind(UR) = 0$. \proofend

We now start to treat general $A$-Fredholm operators. The crucial fact
which allows us to proceed, is that any Fredholm operator is
``regularizable'' over a unital algebra.

\statement{3.8.~Lemma}{Let $B$ be a unital \cstar-algebra and $M$ and
$N$ be Hilbert B-modules. If $T \in \LB(M, N)$ is $B$-Fredholm (but
not necessarily regular), then there exists an integer $n$ and a
regular $B$-Fredholm operator $\Tt$ in $\LB(M\+B^n, N\+B^n)$ such that
$\Tt_{NM} = T$.}

\proofbegin Let $S$ in $\LB(N, M)$ be such that both $I_M - ST$ and
$I_N - TS$ are $B$-finite rank. So, let $\nu = (\nu_i)_i$ and $\mu =
(\mu_i)_i$ be such that $I_M - ST = \Omega_\nu \Omega_\mu^*$. Define
  $$\Tt = \mat{T}{0}{\Omega_\mu^*}{0}
  \quad {\rm and} \quad
  \St = \mat{S}{\Omega_\nu}{0}{0}.$$
  It is easy to see that $\Tt$ and $\St$ are each others
pseudo-inverse hence, in particular, $\Tt$ is regular. The hypothesis
that $B$ have a unit implies that $B^n$ is $B$-finite rank and hence
that
  $$I - \St\Tt = \mat{0}{0}{0}{I_{B^n}}$$ and
  $$I - \Tt\St = \mat{I_N - TS}{-T\Omega_\nu}
  {-\Omega_\mu^*S}{I_{B^n} - \Omega_\mu^* \Omega_\nu}$$
  are $B$-finite rank operators. We conclude that $\Tt$ is Fredholm.
\proofend

Since we shall not assume the algebras we work with to be unital (nor
$\sigma$-unital, as already stressed in the introduction), we will
keep all the applications of the above Lemma to the unitized algebra
$\At$, in the following way.  Given an $A$-Fredholm operator $T$ in
$\LA(M, N)$ consider $M$ and $N$ as Hilbert modules over the unitized
algebra $\At$, as it is done in \cite{\Jensen}, remark 1.1.5.
Obviously $T$ is $\At$-Fredholm as well.  Let therefore $\Tt \in
\LAt(M\+\At^n, N\+\At^n)$ be the operator constructed as in 3.8.
Since $\Tt$ is regular and Fredholm, $\ind(\Tt)$ is well defined as an
element of $K_0(\At)$.

The following result will allow us to return to the realm of
non-unital (meaning non-necessarily-unital) algebras after our brief
encounter with units.

\statement{3.9.~Proposition}{The Fredholm index of $\Tt$ belongs to
$K_0(A)$.}

\proofbegin Denote by $\ep \: \At \arw \C$ the augmentation
homomorphism.  That is $\ep(a + \lambda) = \lambda$ for $\lambda$ in
$\C$ and $a$ in $A$. Since $K_0(A)$ is defined to be the kernel of the
map $\ep_* \: K_0(\At) \arw K_0(\C)$, all we need to do is show that
$\ep_*(\ind(\Tt)) = 0$.

Let $S$, $\St$, $\mu$ and $\nu$ be as in the proof of 3.8.  Note that
$I-\St\Tt = \mat{0}{0}{0}{I_{\At^n}}$, so $\rank(\Ker(\Tt)) = n$. We
then need to show that $\ep_*(\rank(\Ker(\St)))$ is also equal to $n$.
The kernel of $\St$ is the image of the idempotent
  $$I - \Tt\St = \mat{I_N - TS}{-T\Omega_\nu}
  {-\Omega_\mu^*S}{I_{\At^n} - \Omega_\mu^* \Omega_\nu} \eqno
(3.9.1)$$
  which we shall simply denote by $P$. Since $P$ is $\At$-compact,
there are $m$-tuples $\phi = (\phi_1, \ldots, \phi_m)$ and $\psi =
(\psi_1, \ldots, \psi_m)$ of elements of $N\+\At^n$ such that $P =
\Omega_\phi \Omega_\psi^*$. Write each $\phi_i$ as $\phi_i = (\xi_i,
v_i)$ and $\psi_i = (\eta_i, w_i)$.

Replacing, if necessary, $\phi_i$ by $P(\phi_i)$, we can assume that
$P\Omega_\phi = \Omega_\phi$ and therefore that $Q = \Omega_\psi^*
\Omega_\phi$ is an idempotent operator on $\At^n$ whose image is
isomorphic, as $\At$-modules, to the image of $P$. As an $m\times m$
matrix, $Q$ is given by $Q = \bigl(\<\psi_i,\phi_j\>\bigr)_{i, j}$.
Our goal is then to show that the trace of the complex idempotent
matrix $\ep(Q) = \bigl(\ep\bigl(\<\psi_i,\phi_j\>\bigr)\bigr)_{i, j}$
equals $n$. That trace is given by
  $$\sum_{i=1}^m \ep\bigl(\<\psi_i,\phi_i\>\bigr)
  = \sum_{i=1}^m \ep\bigl(\<w_i, v_i\>\bigr)
  = \ep\biggl(\sum_{i=1}^m \sum_{r=1}^n \<w_i, e_r\> \<e_r,
v_i\>\biggr), $$
  where $\{e_i, \ldots, e_n\}$ is the canonical basis of $\At^n$. The
above then equals
  $$\ep\biggl(\sum_{i=1}^m \sum_{r=1}^n \<e_r, v_i\<w_i,
e_r\>\>\biggr)
  = \ep\biggl(\sum_{r=1}^n \<(0, e_r), P(0, e_r)\>\biggr).$$
  Using the definition of $P$ in 3.9.1, the term $\sum_{r=1}^n \<(0,
e_r), P(0, e_r)\>$ can be expressed as
  $$\sum_{r=1}^n \<e_r, (I - \Omega_\mu^* \Omega_\nu)e_r\>
  = n - \sum_{r=1}^n \<\mu_r,\nu_r\>.$$
  The last term above clearly maps to $n$ under $\ep$ so the proof is
complete. \proofend

The statement of 3.9 is meant to refer to the specific construction of
$\Tt$ obtained is 3.8. But note that any regular Fredholm operator in
$\LAt(M\+\At^n, N\+\At^n)$, which has $T$ in the upper left corner,
will differ from the $\Tt$ above, by an $\At$-compact operator.
Therefore its index will coincide with that of $\Tt$ by 3.7, and so
will be in $K_0(A)$ as well.

\statement{3.10.~Definition}{If $T$ is an $A$-Fredholm operator in
$\LA(M, N)$, then the Fredholm index of $T$, denoted $\ind(T)$, is
defined to be the index of the regular Fredholm operator $\Tt$
constructed in 3.8.}

Clearly, if $T$ is already regular, we can take $n=0$ in 3.8 so that
the above definition extends the one given in 3.4.  Elementary
properties of the Fredholm index are collected in the next
proposition.

\statement{3.11.~Proposition}{If $T \in \LA(M, N)$ is a Fredholm
operator, then
  \medskip \item{(i)} $\ind(T^*) = - \ind(T)$
  \medskip \item{(ii)} If $U$ in $\LA(X, M)$ and $V$ in $\LA(Y, N)$
are invertible, then $\ind(VTU) = \ind(T)$.
  \medskip \item{(iii)} If $T'\in\LA(M, N)$ is such that $T' - T$ is
in $\KA(M, N)$, then $T'$ is also Fredholm and $\ind(T') = \ind(T)$.
  \medskip \item{(iv)} If $S \in \LA(N, M)$ is such that $I_M - ST$ is
in $\KA(M)$, then $\ind(S) = -\ind(T)$.
  \medskip \item{(v)} If $T_1 \in \LA(M_1, N_1)$ is Fredholm, then
$\ind(T\+T_1) = \ind(T) + \ind(T_1)$.}

\proofbegin Left to the reader. \proofend

Let us now briefly tackle the question of invariance of the Fredholm
index under small perturbations. The proof, which we omit, is similar
to the classical one.

\statement{3.12.~Proposition}{Let $T$ in $\LA(M, N)$ be Fredholm. Then
there is a positive real number $\ep$ such that any $T'$ satisfying
$\|T'-T\| \less \ep$ is also Fredholm with $\ind(T') = \ind(T)$. In
fact $\ep$ can be taken to be $\ep = \|S\|^{-1}$ for any $S$ in
$\LA(N, M)$ such that $I_M - ST$ is in $\KA(M)$.}

It is equally easy to prove:

\statement{3.13.~Proposition}{If $T_1 \in \LA(M, N)$ and $T_2 \in
\LA(N, P)$ are Fredholm operators, then $T_2T_1$ is Fredholm and
$\ind(T_2T_1) = \ind(T_2) + \ind(T_1)$.}

It is part of our goal to find an alternate definition of $K_0(A)$ in
terms of Fredholm operators. For this reason it is important to have a
sufficiently large collection of such operators. Specifically we will
need to exhibit Fredholm operators with an arbitrary element of
$K_0(A)$ as the index.

This is quite easy if $A$ has a unit: given an arbitrary element
$[M]_0 - [N]_0$ in $K_0(A)$, where $M$ and $N$ are finitely generated
projective $A$-modules, consider the zero operator ${\bf 0}\: M \arw
N$.  If both $M$ and $N$ are given the Hilbert module structure they
possess, being isomorphic to direct summands of $A^n$, it follows that
{\bf 0} is a Fredholm operator, its index obviously being $[M]_0 -
[N]_0$.  The situation is slightly more complicated in the non-unital
case, and it is treated in our next result.

\statement{3.14.~Proposition}{For any $\alpha$ in $K_0(A)$ there is an
$A$-Fredholm operator $T$ such that $\ind(T) = \alpha$.}

\proofbegin Regarding $K_0(A)$ as the kernel of the augmentation map
$\ep_*$, in $K_0(\At)$, write $\alpha = [p]_0 - [q]_0$, where $p$ and
$q$ are self-adjoint idempotent $n\times n$ matrices over $\At$ such
that $\ep_*([p]_0 - [q]_0) = 0$. It follows that $\ep(p)$ and $\ep(q)$
are similar complex matrices. Hence, after performing a conjugation
of, say $q$, by a complex unitary matrix, we may assume that $\ep(p)$
and $\ep(q)$ are in fact equal. Therefore $p-q$ is in $M_n(A)$.

We now claim that the operator $T\:pA^n \arw qA^n$ given by $T(v) =
qv$ is Fredholm and that its index is $[p]_0 - [q]_0$.

Let $S\:qA^n \arw pA^n$ be given by $S(v) = pv$.  Denote by $p_i$ the
$i^{\rm th}$ column of $p$, viewed as a $n\times 1$ matrix, and let
$(u_\lambda)_\lambda$ be an approximate identity for $A$. Define
  $\xi_i = p(p-q)p_i$ and
  $\eta_i^\lambda = p_i u_\lambda$
  for $1 \leq i \leq n$.
  So, $\xi_i$ and $\eta_i^\lambda$ are elements of $pA^n$.
  Let $\xi = (\xi_i)_i$ and $\eta^\lambda = (\eta_i^\lambda)_i$. For
$\zeta$ in $pA^n$ we have
  $$\Omega_\xi\Omega_{\eta^\lambda}^*(\zeta)
  = \sum_{i=1}^n \xi_i \<\eta_i^\lambda, \zeta\>
  = \sum_{i=1}^n p(p-q)p_iu_\lambda p_i^* \zeta$$ which converges,
uniformly in $\|\zeta\| \leq 1$, as $\lambda \arw \infty$, to
  $$\sum_{i=1}^n p(p-q)p_i p_i^* \zeta
  = p(p-q)p\zeta
  = \zeta - pq\zeta
  = (I - ST)\zeta.$$
  In a similar fashion we can show that $I-TS$ is $A$-compact as well.
In order to compute the index of $T$, consider the operators
  $\Tt$ in $\LAt(pA^n\+\At^n, qA^n\+\At^n)$ and
  $\St$ in $\LAt(qA^n\+\At^n, pA^n\+\At^n)$ given by
  $$\Tt = \mat{qp} {q(I-p)} {(I-q)p} {(I-q)(I-p)}
  \quad {\rm and} \quad
  \St = \mat{pq} {p(I-q)} {(I-p)q} {(I-p)(I-q)}.$$
  Direct computation shows that
  $$\St\Tt = \mat{I} {0} {0} {I-p}
  \quad {\rm and} \quad
  \Tt\St = \mat{I} {0} {0} {I-q}, $$
  from which it follows that $\St$ is a pseudo-inverse for $\Tt$ and
hence that $\Tt$ is a regular $\At$-Fredholm operator.  By definition
we have
  $$\ind(T) = \ind(\Tt) = \rank(I - \St\Tt) - \rank(I - \Tt\St)
  = [p]_0 - [q]_0.\eqno\square$$ \medbreak

Our last result showed that any $K_0$ element is the index of some
Fredholm operator. We would now like to discuss the question of when
two Fredholm operators have the same index. As a first step, we
classify operators with index zero.

\statement{3.15.~Lemma}{Let $T \in \LA(M, N)$ be a Fredholm operator
with $\ind(T) = 0$. Then there exists an integer $n$ such that
$T\+I_{A^n} \: M\+A^n \arw N\+A^n$ is an $A$-compact perturbation of
an invertible operator.}

\proofbegin Let $\Tt$ in $\LAt(M\+\At^n, N\+\At^n)$ be constructed as
in 3.8 with $\At$ playing the role of $B$. If $\St$ is also as in 3.8
we have $I-\St\Tt = \mat{0}{0}{0}{I_{\At^n}}$. The hypothesis that
$\ind(T) = 0$ then says that $\rank(I - \Tt\St) = \rank(I - \St\Tt)$.
But, as we see above, the rank of $I-\St\Tt$ is the same as the rank
of the free $\At$-module $\At^n$. So, we have that $\Im(I-\Tt\St)$ is
stably isomorphic, as $\At$-modules, to $\At^n$. By increasing $n$, if
necessary, we may thus assume that $\Im(I-\Tt\St)$ is, in fact,
isomorphic to $\At^n$.

This said, we may find a generating set $\{x_i\}_{i=1}^n$ for
$\Im(I-\Tt\St)$ which is orthonormal in the sense that $\<x_i, x_j\> =
\delta_{ij}$. Let each $x_i$ be given by $x_i = (\xi_i, v_i)$. Each
$v_j$ is in $\At^n$ and hence we may write $v_j = (v_{ij})_i$ with
$v_{ij} \in \At$. The fact that the $x_i$ form an orthonormal set
translates to
  $$\<\xi_i,\xi_j\> + \sum_{k=1}^n v_{ki}^*v_{kj} = \delta_{ij}.$$

Recalling that $\ep \: \At \arw \C$ denotes the augmentation
homomorphism, observe that the complex matrix $u =
\bigl(u_{ij}\bigr)_{i, j}$, where $u_{ij} = \ep(v_{ij})$, is unitary.
If we now define $x'_j = \sum_{i=1}^n u_{ji}^* x_i$ and write $x'_j =
(\xi'_j, v'_j)$ with $v'_j = (v'_{ij})_i$ one can show that
$\ep(v'_{ij}) = \delta_{ij}$. In other words, we may assume, without
loss of generality, that $\ep(v_{ij}) = \delta_{ij}$. Otherwise,
replace each $x_i$ by $x'_i$.

Let $\xi = (\xi_i)_i$ and $v = (v_i)_i$ so that $\Omega_\xi \+
\Omega_v$ is the isomorphism from $\At^n$ onto $\Im(I-\Tt\St)$
mentioned above. The operator $U$ in $\LAt(M\+\At^n, N\+\At^n)$ given
by
  $$U = \mat{T}{\Omega_\xi}{\Omega_\mu^*}{\Omega_v}$$ is therefore
invertible (please note that $\mu$ is as in the proof of 3.8).

Note that
  $$\overline{(M\+\At^n)\cdot A} = M\+A^n$$
  and similarly
  $$\overline{(N\+\At^n)\cdot A} = N\+A^n.$$
  So, $U$ gives, by restriction, an invertible operator from $M\+A^n$
to $N\+A^n$ and, denoting the latter by $U$, from now on, we have
  $$U - \mat{T}{0}{0}{I_{A^n}}
  = \mat{0}{\Omega_\xi}{\Omega_\mu^*}{\Omega_v - I_{A^n}}.$$
  The matrix on the right hand side represents an $A$-compact
operator: the crucial point being that $\Omega_v - I_{A^n}$ is the
operator on $A^n$ given by multiplication by the matrix $\bigl(v_{ij}
- \delta_{ij}\bigr)_{i, j}$ which is in $M_n(A)$ since $\ep(v_{ij} -
\delta_{ij})$ was seen to be zero (see \cite{\KI}, Lemma 2.4).
\proofend

The following characterization of when two Fredholm operators have the
same index is an immediate corollary of our last Lemma.

\statement{3.16.~Proposition}{If $T_i$ in $\LA(M_i, N_i)$ for $i=1,
2$, are Fredholm operators such that $\ind(T_1) = \ind(T_2)$, then for
some integer $n$, the operator
  $$T_1 \+ T_2^* \+ I_{A^n}\ :\ M_1 \+ N_2 \+ A^n \arw N_1 \+ M_2 \+
A^n$$
  is an $A$-compact perturbation of an invertible operator.}

Using the machinery developed so far, we may provide an alternate
definition of the $K$-theory group $K_0(A)$. Choose, once and for all,
a cardinal number $\omega$ which is bigger than the cardinality of
$A^n$ for every integer $n$.  We remark that the role of $\omega$ is
merely to avoid set theoretical problems arising from the careless
reference to the set of {\it all\/} $A$-Fredholm operator.  Any choice
of $\omega$, as long as it is sufficiently large, will result in the
same conclusions.

Denote by $F_0(A)$ the set of all $A$-Fredholm operators whose domain
and codomain are Hilbert modules of cardinality no larger than
$\omega$ (actually we should require these Hilbert modules to have a
subset of $\omega$ as their carrier set). Declare two elements $T_1$
and $T_2$ of $F_0(A)$ equivalent, if there is an integer $n$ such that
$T_1 \+ T_2^* \+ I_{A^n}$ is an $A$-compact perturbation of an
invertible operator.

The quotient $F(A)$ of $F_0(A)$ by the above equivalence relation is
obviously a group with the operation of direct sum of Fredholm
operators. The inverse of the class of $T$ being given by that of
$T^*$ by 3.11 and 3.15.

\statement{3.17.~Corollary}{The Fredholm index map, viewed as a map
  $$\ind \: F(A) \arw K_0(A), $$
  is a group isomorphism.}

\proofbegin Follows immediately from 3.14 and 3.16. \proofend

We should remark that 3.17 is a generalization of the fact that
$\KK(\C, A)$ is isomorphic to $K_0(A)$. The new aspect being that no
separability is involved.  This is one of the crucial steps in
achieving our main result as we shall see shortly.

\beginsection 4. Preliminaries on Hilbert Bimodules

We would now like to set the present section aside in order to present
a few relevant aspects of the theory of Hilbert bimodules which will
be important for our discussion of Morita equivalence. We adopt the
definition of Hilbert bimodules given in \cite{\BMS}, 1.8. Namely, if
$A$ and $B$ are \cstar-algebras, a Hilbert $A$-$B$-bimodule is a
complex vector space $X$ which is a left Hilbert $A$-module as well as
a right Hilbert $B$-module, and such that the $A$-valued inner product
  $$(\cdot|\cdot)\:X\times X\arw A$$
  and the $B$-valued inner product
  $$\<\cdot\thinspace,\cdot\>\:X\times X\arw B$$
  are related by the identity
  $$(\xi|\eta)\mu = \xi\<\eta,\mu\>, \quad \xi,\eta,\mu\in X.$$

Some authors prefer to use the notation $\<\cdot\thinspace,\cdot\>_A$
and $\<\cdot\thinspace,\cdot\>_B$ for these inner-products but we
believe the notation indicated above makes some formulas much more
readable. In particular, it is implicit that any inner-product denoted
by $\<\cdot\thinspace,\cdot\>$ will be linear in the second variable
while those denoted by $(\cdot|\cdot)$ are linear in the first
variable. We should nevertheless remark that the differentiated
notation is not meant to imply any asymmetry in the structure of
bimodules. With the obvious interchange of left and right, any result
that holds on the ``left'' will also hold on the ``right'' and
vice-versa.

As mentioned in \cite{\BMS}, Hilbert $A$-$B$-bimodules are nothing but
Rieffel's imprimitivity bimodules (see \cite{\RI}, 6.10) for which it
is {\it not\/} assumed that the range of the inner-products generate
the coefficient algebras.

The closed span of the set
  $\{\<x, y\>\:x, y\in X\}$,
  which we denote by $\<X, X\>$, is a two sided ideal in $B$ and
similarly, the closed span of
  $\{(x|y)\:x, y\in X\}$ is the ideal $(X|X)$ of $A$.

\statement{4.1.~Definition}{The Hilbert $A$-$B$-bimodule $X$ is said
to be {\it left-full} (resp.~{\it right-full}) if $(X|X)$ coincides
with $A$ (resp.~if $\<X, X\>$ coincides with $B$).}

Using the terminology just introduced, Rieffel's imprimitivity
bimodules are precisely the Hilbert-bimodules that are simultaneously
left-full and right-full.

Throughout this section we shall consider fixed two \cstar-algebras
$A$ and $B$ as well as a Hilbert $A$-$B$-bimodule $X$. As before, $M$
and $N$ will denote Hilbert $A$-modules.  If $M$ is a Hilbert
$A$-module (we remind the reader of our convention according to which
{\it module\/} without further adjectives, means {\it right
module\/}), then the tensor product module $M\*_A X$ has a natural
$B$-valued (possibly degenerated) inner-product specified by
  $$\<\xi_1\*x_1,\xi_2\*x_2\> = \<x_1,\<\xi_1,\xi_2\>x_2\>, \quad
\xi_1,\xi_2\in M,\ x_1, x_2\in X.$$
  After moding out the elements of norm zero and completing, we are
left with a Hilbert $B$-module which we also denote, for simplicity,
by $M\*_A X$. See \cite{\Jensen}, 1.2.3 for details, but please
observe that the notation used there is not the same as the one just
described.  It should also be observed that one does not need the
$A$-valued inner-product on $X$ in order to perform this construction.
It is enough that $X$ be a left $A$-module in such a way that the
representation of $A$, as left multiplication operators on $X$, be a
$*$-homomorphism.

If $T$ is in $\LA(M, N)$, we denote by $T\*I_X$ the linear
transformation
  $$T\*I_X \ :\ M\*_A X \arw N\*_A X$$
  given by $T\*I_X(\xi, x) = T(\xi)\*x$ for $\xi$ in $M$ and $x$ in
$X$.  A slight modification of \cite{\Jensen}, 1.2.3 shows that
$T\*I_X$ is in $\LB(M\*_AX, N\*_AX)$ and that $\|T\*I_X\| \leq \|T\|$.

Let us now present one of our most important technical results.
Although quite a simple fact, with an equally simple proof, it is a
crucial ingredient in this work.  Compare \cite{\Jensen}, 1.2.8.

\statement{4.2.~Theorem}{If the Hilbert $A$-$B$-bimodule $X$ is
left-full and if $T$ is in $\KA(M, N)$, then $T\*I_X$ is $B$-compact.}

\proofbegin It obviously suffices to prove the statement in case
$T=\Omega_\nu \Omega_\mu^*$ with $\mu$ in $M$ and $\nu$ in $N$. Given
$\xi\*x$ in $M\*_A X$ we have
  $$T\*I_X(\xi\*x) = \nu\<\mu,\xi\>\*x = \nu\*\<\mu,\xi\>x.$$ Observe
that, since $X$ is left-full and also by \cite{\Jensen}, Lemma 1.1.4,
there is no harm in assuming that $\nu = \nu_1(y|z)$ for some $y, z$
in $X$ and $\nu_1$ in $N$. So
  $$T\*I_X(\xi\*x) = \nu_1\*(y|z)\<\mu,\xi\>x
  = \nu_1\*y\<z,\<\mu,\xi\>x\>$$
  $$= (\nu_1\*y)\<\mu\*z,\xi\*x\> = \Omega_{\nu_1\*y}
\Omega_{\mu\*z}^*(\xi\*x).$$
  This concludes the proof. \proofend

One of the main uses we shall have for this result is recorded in:

\statement{4.3.~Corollary}{If $X$ is left-full and if $T\in\LA(M, N)$
is an $A$-Fredholm operator, then $T\*I_X$ is $B$-Fredholm.}

\proofbegin If $S\in\LA(N, N)$ is such that $I_M-ST$ is in $\KA(M)$
and $I_N-TS$ is in $\KA(N)$, then
  $$I_{M\*_A X} - (S\*I_X)(T\*I_X) = (I_M -ST)\*I_X$$
  which is a $B$-compact operator by 4.2. The same reasoning applies
to
  $I_{M\*_A X} - (T\*I_X)(S\*I_X)$. \proofend

At this point, the reader may have already anticipated our strategy of
using a bimodule to create a homomorphism on $K_0$-groups: given an
element $\alpha$ in $K_0(A)$, we may find, using 3.14, an $A$-Fredholm
operator $T$ whose index is $\alpha$. The index, in $K_0(B)$, of the
$B$-Fredholm operator $T\*I_X$ is the image of $\alpha$ under the
homomorphism we have in mind. In order to make this picture work, we
need to tackle the question of well definedness, which we now do.

\statement{4.4.~Proposition}{If $T_1, T_2\in\LA(M, N)$ are
$A$-Fredholm operators such that $\ind(T_1) = \ind(T_2)$, then
$\ind(T_1\*I_X) = \ind(T_2\*I_X) $.}

\proofbegin According to 3.16 there is an integer $n$ such that
  $$T_1 \+ T_2^* \+ I_{A^n} \ :\ M_1 \+ N_2 \+ A^n \arw N_1 \+ M_2 \+
A^n$$
  is an $A$-compact perturbation of an invertible operator. By 4.2 and
by the fact that the tensor product distributes over direct sums, we
have that
  $(T_1\*I_X) \+ (T_2^*\*I_X) \+ (I_{A^n}\*I_X)$
  is a $B$-compact perturbation of an invertible operator. By 3.11 its
index is therefore zero. On the other hand, also by 3.11 we have
  $$\ind(T_1\*I_X) + \ind(T_2^*\*I_X) + \ind(I_{A^n}\*I_X) = 0$$
  which says that $\ind(T_1\*I_X) = \ind(T_2\*I_X)$. \proofend

A important ingredient for the functoriality properties of left-full
Hilbert bimodules is the notion of tensor product of bimodules. In
order to avoid endless calculations that arise in an abstract
treatment of tensor products, we shall provide an alternative picture
for bimodules, as concrete operators between Hilbert spaces, in which
the coefficient algebras are represented.  The notion of
representation of bimodules is described next. Compare \cite{\BMS},
Definition 2.1.

\statement{4.5.~Definition}{A {\it representation} of the Hilbert
$A$-$B$-bimodule $X$ consists of the following data:
  \medskip\item{(i)} a representation $\pi_A$ of $A$ on a Hilbert
space $H_A$,
  \medskip\item{(ii)} a representation $\pi_B$ of $B$ on a Hilbert
space $H_B$ and
  \medskip\item{(iii)} a bounded linear transformation $\pi_X$ from
$X$ into the Banach space $\B(H_B, H_A)$ of all bounded linear
operators from $H_B$ to $H_A$.
  \medskip Furthermore it is required that, for $a\in A$, $b\in B$,
and $x, x_1, x_2\in X$,
  \medskip\item{(a)} $\pi_X(ax)=\pi_A(a)\pi_X(x)$
  \medskip\item{(b)} $\pi_X(xb)=\pi_X(x)\pi_B(b)$
  \medskip\item{(c)} $\pi_A\bigl((x_1|x_2)\bigr) =
\pi_X(x_1)\pi_X(x_2)^*$
  \medskip\item{(d)} $\pi_B\bigl(\<x_1, x_2\>\bigr) =
\pi_X(x_1)^*\pi_X(x_2)$.}

At this point it is necessary to remark that for $x$ in $X$, one has
that $\|(x|x)\| = \|\<x, x\>\|$ (see \cite{\BMS}, Remark 1.9).  So,
when we speak of $\|x\|$, we mean the square root of that common
value. In particular, this is the norm we have in mind when we
require, in (iii), that $\pi_X$ be a bounded map on $X$.

\statement{4.6.~Proposition}{Let $(\pi_A,\pi_B,\pi_X)$ be any
representation of $X$. If either $\pi_A$ or $\pi_B$ are faithful, then
$\pi_X$ is isometric.}

\proofbegin Let $x\in X$. We have
  $\|\pi_X(x)\|^2 = \|\pi_X(x)\pi_X(x)^*\| =
\|\pi_A\bigl((x|x)\bigr)\|$.
  So, assuming that $\pi_A$ is faithful, we have
  $\|\pi_X(x)\|^2 =\|(x|x)\| = \|x\|^2$.
  A similar argument applies if $\pi_B$ is assumed to be faithful,
instead. \proofend

Given a representation $\pi_B$ of $B$, it is natural to ask whether or
not $\pi_B$ is part of the data forming a representation of $X$. To
answer this question we need to bring in the conjugate module and the
linking algebra.  The conjugate module of $X$ is the bimodule one
obtains by reversing its structure so as to produce a Hilbert
$B$-$A$-bimodule as explained in \cite{\RI}, 6.17, or \cite{\BMS},
1.4.  We shall denote the conjugate module by $X^*$ (although
$\tilde{X}$ is used in \cite{\RI}).  The linking algebra of $X$,
introduced in \cite{\BGR}, 1.1 in the special case when $X$ is both
left and right-full, and in \cite{\BMS}, 2.2 in general, is the
\cstar-algebra
  $$ L = \mat{A} {X} {X^*} {B}$$
  equipped with the multiplication
  $$\mat{a_1} {x_1} {y_1^*} {b_1} \mat{a_2} {x_2} {y_2^*} {b_2}
  = \mat{a_1 a_2 + (x_1|y_2)} {a_1 x_2 + x_1 b_2}
    {y_1^* a_2 + b_1 y_2^*} {\<y_1, x_2\> + b_1 b_2}$$
  and involution
  $$\mat{a} {x} {y^*} {b} = \mat{a^*} {y} {x^*} {b^*}$$
  for $a, a_1, a_2 \in A$, $b, b_1, b_2 \in B$ and $x, x_1, x_2, y,
y_1, y_2\in X$.  Here, $x^*$ denotes the element $x$ of $X$ when it is
viewed in $X^*$.

\statement{4.7.~Proposition}{Let $\pi_B$ be a non-degenerated
representation of $B$ on the Hilbert space $H_B$. Then there exists a
Hilbert space $H_A$, a non-degenerated representation $\pi_A$ of $A$
on $H_A$ and a bounded linear map $\pi_X\: X\arw \B(H_B, H_A)$ which,
when put together, form a representation of $X$.}

\proofbegin Let $L$ be the linking algebra of $X$. Thus $\pi_B$
becomes a representation of a subalgebra of $L$, namely $B$. Let $\pi$
be a representation of $L$ on a Hilbert space $H$ which
  contains a copy of $H_B$,
  such that $H_B$ is invariant under the operators $\pi(b)$, for $b$
in $B$,
  and, finally, such that $\pi(b)|_{H_B}=\pi_B(b)$.
  The existence of $\pi$, in the case that $\pi_B$ is cyclic, follows
{}from the Theorem on extension of states and the GNS construction. In
the general case, it follows from the fact that any representation is
a direct sum of cyclic representations.

Viewing $X$ as the subset of $L$ formed by the matrices with the only
non-zero entry lying in the upper right hand corner, let $H_A =
\overline{\pi(X) H_B}$. Here, and in the sequel, products of sets, as
in ``$\pi(X) H_B$'', will always mean the linear span of the set of
individual products.

Since $X$ is a left $A$-module, it is clear that $H_A$ is invariant
under $\pi(A)$. Denote by $\pi_A$ the sub-representation of $A$ on
$H_A$ given in this way. Observe that $\overline{AX} = X$, by
\cite{\Jensen}, 1.1.4, so we conclude that $\pi_A$ is non-degenerated.

For each $x$ in $X$ we have, by definition, that $\pi(x)H_B \subseteq
H_A$. By totally different reasons we also have that $\pi_X(x)^*
H_A\subseteq H_B$. In fact, let $\xi\in H_A$. Without loss of
generality we may assume that $\xi=\pi(y)\eta$ where $y$ is in $X$ and
$\eta\in H_B$. We then have
  $$\pi(x)^*\xi =\pi(x)^*\pi(y)\eta = \pi\bigl(\<x, y\>\bigr)\eta \in
H_B.$$
  For each $x$ in $X$ let $\pi_X(x)$ denote the element of $\B(H_B,
H_A)$ obtained by restriction of $\pi(x)$. The properties of
Definition 4.5 may now be easily checked. If the reader does decide to
do so, we suggest the formal definition $\pi_X(x) := i_A^*\pi(x)i_B$
where $i_A$ and $i_B$ are the inclusion operators from $H_A$ and $H_B$
into $H$. This has the advantage of taking care of the subtle issue of
reducing the size of the co-domain of an operator. \proofend

\statement{4.8.~Proposition}{There exists a representation $(\pi_A,
\pi_B, \pi_X)$ such that both $\pi_A$ and $\pi_B$ are faithful (and
hence $\pi_X$ is isometric by 4.6).}

\proofbegin Let $\pi$ be a faithful representation of the linking
algebra $L$ on the Hilbert space $H$. Define $H_A = \overline{\pi(A)
H}$ and $H_B = \overline{\pi(B) H}$ and let $\pi_A$ and $\pi_B$ be the
corresponding sub-representations of $A$ and $B$ on $H_A$ and $H_B$,
respectively. Since $\overline{AX} = X$ and $\overline{BX^*} = X^*$ it
is clear that $\pi(X)H \subseteq H_A$ and that $\pi(X)^*H\subseteq
H_B$.  If, for each $x$ in $X$, we denote by $\pi_X(x)$ the element of
$\B(H_B, H_A)$ given by restriction of $\pi(x)$, the proof can be
completed as in 4.7. \proofend

\statement{4.9.~Lemma}{The set
  $\bigl\{\sum_{i=1}^n(x_i|x_i)\: n\in\N,\ x_i\in X\bigr\}$
  is dense in the positive cone of $(A|A)$.}

\proofbegin For $a=\sum_{i=1}^n(x_i|y_i)$. We have
  $$a^*a = \sum_{i, j}(y_i|x_i)(x_j|y_j)
  = \sum_{i, j}\bigl((y_i|x_i)x_j|y_j\bigr)
  = \sum_{i, j}\bigl(y_i\<x_i, x_j\>|y_j\bigr).$$
  The matrix $\bigl(\<x_i, x_j\>\bigr)_{i, j} \in M_n(B)$ is a
positive matrix as observed in \cite{\Jensen}, 1.2.4.  So, there is
$\bigl(b_{ij}\bigr)_{i, j}$ in $M_n(B)$ such that
  $$\<x_i, x_j\> = \sum_{k=1}^n b_{ik} b_{jk}^*.$$
  We then have
  $$a^*a = \sum_{i, j, k}(y_ib_{ik}|y_jb_{jk}).$$
  If we now define $z_k = \sum_{i=1}^n y_ib_{ik}$ we have
  $$a^*a = \sum_{k=1}^n(z_k|z_k).$$
  Since the set of elements $a^*a$, with $a$ as above, is clearly
dense in the positive cone of $(A|A)$, the proof is complete.
\proofend

\statement{4.10.~Proposition}{Let $(\pi_A,\pi_B,\pi_X)$ be a
representation of $X$. If $\pi_B$ is faithful then $\pi_A$ is faithful
on $(A|A)$.}

\proofbegin Using 4.9, it is enough to show that
  $\|\pi_A\bigl(\sum_{i=1}^n (x_i|x_i)\bigr)\|
  = \|\sum_{i=1}^n (x_i|x_i)\|$.
  Thus, let $h = \sum_{i=1}^n (x_i|x_i)$ and observe that
  $$\|\pi_A(h)\|
  = \|\sum_{i=1}^n \pi_X(x_i)\pi_X(x_i)^*\|
  = \left\Vert
    \bigl(\pi_X(x_1), \ldots, \pi_X(x_n)\bigr)
    \pmatrix{\pi_X(x_1)^* \cr \vdots \cr \pi_X(x_n)^*}
  \right\Vert$$
  where, in the last term above, we mean the product of a row matrix
by a column matrix. Since the identity $\|T^*T\| = \|TT^*\|$ holds for
general operators, the above equals
  $$\left\Vert
    \pmatrix{\pi_X(x_1)^* \cr \vdots \cr \pi_X(x_n)^*}
    \bigl(\pi_X(x_1), \ldots, \pi_X(x_n)\bigr)
  \right\Vert
  = \|\bigl(\pi_X(x_i)^*\pi_X(x_j)\bigr)_{i, j}\|$$
  $$ = \|\bigl(\pi_B\bigl(\<x_i, x_j\>\bigr)\bigr)_{i, j}\|
  = \|\bigl(\<x_i, x_j\>\bigr)_{i, j}\|.$$
  The exact same computations done so far can obviously be repeated,
in reverse order, for a representation $(\rho_A, \rho_B, \rho_X)$ in
which all components are isometric, as for example, the representation
provided by 4.8.  This shows that the last term above equals
$\|\rho_A(h)\| = \|h\|$. So, $\|\pi_A(h)\| = \|h\|$.  \proofend

{}From this point on, and until the end of this section, we shall
consider fixed another \cstar-algebra, denoted $C$, and a Hilbert
$B$-$C$-bimodule $Y$. Our goal is to make sense of $X\*_BY$ as a
Hilbert $A$-$C$-bimodule. So, for the time being, let us denote by
$X\*_BY$, the algebraic tensor product of $X$ and $Y$ over $B$ which
provides us with an $A$-$C$-bimodule.

\statement{4.11.~Definition}{Let
  $(\cdot|\cdot)$ and $\<\cdot\thinspace,\cdot\>$ be the sesqui-linear
forms on $X\*_BY$ (the first one being linear in the first variable
and vice-versa) specified by
  $$(x_1\*y_1|x_2\*y_2) = (x_1(y_1|y_2)|x_2)$$
  and
  $$\<x_1\*y_1, x_2\*y_2\> = \<y_1,\<x_1, x_2\>y_2\>.$$}

The only steps that are not entirely trivial in checking that this can
be made into a Hilbert $B$-$C$-bimodule are that
  \medskip\item{(a)} both sesqui-linear forms above are positive and
  \medskip\item{(b)} $\|(z|z)\| = \|\<z, z\>\|$ for all $z$ in
$X\*_BY$.
  \medskip In order to check these, fix a faithful non-degenerated
representation $\pi_B$ of $B$ in some Hilbert space $H_B$. Using 4.7
we may find a representation $\pi_A$ on a space $H_A$ as well as a
representation $\pi_X$ of $X$, by bounded operators from $H_B$ to
$H_A$, satisfying the axioms described in Definition 4.5.

Using the symmetric version of 4.7 we can also find a Hilbert space
$H_C$ as well as $\pi_C$ and $\pi_Y$ satisfying the conditions of 4.5.
Note that, by 4.6, both $\pi_X$ and $\pi_Y$ are isometric. By 4.10 it
follows that $\pi_A$ is faithful on $(A|A)$ and that $\pi_C$ is
faithful on $\<C, C\>$. Consider the map
  $$\rho \: X\*_BY \arw \B(H_C, H_A)$$
  given by $\rho(x\*y) = \pi_X(x)\pi_Y(y)$ (meaning composition of
operators). Observe that, for $x_1$, $x_2$ in $X$ and $y_1$, $y_2$ in
$Y$, we have
  $$\pi_A\bigl((x_1\*y_1|x_2\*y_2)\bigr)
  = \pi_X(x_1(y_1|y_2)) \pi_X(x_2)^*$$
  $$ = \pi_X(x_1) \pi_Y(y_1) \pi_Y(y_2)^* \pi_X(x_2)^*
  = \rho(x_1\*y_1) \rho(x_2\*y_2)^*.$$
  Similarly we have
  $\pi_C\bigl(\<x_1\*y_1, x_2\*y_2\>\bigr)
  = \rho(x_1\*y_1)^* \rho(x_2\*y_2)$.
  If $z = \sum_{i=1}^n x_i\*y_i$ is an arbitrary element of $X\*_BY$
we then have that
  $$\pi_A\bigl((z|z)\bigr)
  = \left(\sum\rho(x_i\*y_i)\right) \left(\sum\rho(x_i \*y_i)
\right)^*$$
  which is clearly a positive element in $\B(H_A)$. Therefore, since
$\pi_A$ is faithful on $(X|X)$, we conclude that $(z|z)$ is positive.
Similarly it can be shown that $\<z, z\>$ is positive as well. This
proves (a) above. With respect to (c) we have
  $$\|(z|z)\| = \|\pi_A\bigl((z|z)\bigr)\|
  = \left\Vert \left(\sum\rho(x_i\*y_i)\right)
               \left(\sum\rho(x_i\*y_i)\right)^* \right\Vert$$
  $$ = \left\Vert \left(\sum\rho(x_i\*y_i)\right)^*
               \left(\sum\rho(x_i\*y_i)\right) \right\Vert
  = \|\pi_C\bigl(\<z, z\>\bigr)\| = \|\<z, z\>\|.$$

This said, we may define unambiguously, a semi-norm in $X\*_BY$ by
  $\|z\| = \|(z|z)\| = \|\<z, z\>\|$.  After moding out by the
elements of norm zero and completing, we are left with a Hilbert
$A$-$C$-bimodule, which we denote, for simplicity, by $X\*_BY$, as
well. Observe that the procedure of moding out the null elements and
completing is equivalent to considering the closure of $\rho(X\*_BY)$
in $\B(H_C, H_A)$, which, incidentally is the same as
$\overline{\pi_X(X) \pi_Y(Y)}$.

In light of 4.5, it is easy to see that the $A$-$C$-bimodule structure
is reproduced as the usual composition of operators.  Furthermore, the
$A$-valued inner product becomes $(T|S) = TS^*$ while $\<T, S\> =
T^*S$, as long as the appropriate identifications are made.

\statement{4.12.~Proposition}{If both $X$ and $Y$ are left-full then
so is $X\*_BY$.}

\proofbegin By \cite{\Jensen}, 1.1.4, it follows that $\overline{XB} =
X$. But, since we are assuming that $(Y|Y)=B$ we get
$\overline{X(Y|Y)} = X$. It follows that
  $$(X\*_BY|X\*_BY) = \overline{(X(Y|Y)|X)} = (X|X) = A.\eqno\square$$
\medbreak

Recall that $X^*$ denotes the conjugate module of $X$. Clearly $X^*$
is left-full (resp.~right-full) if and only if $X$ is right-full
(resp.~left-full).

\statement{4.13~Proposition}{The tensor product Hilbert $A$-$A$
bimodule $X\*_BX^*$ is isomorphic to $(A|A)$ (once $(A|A)$ is given
its obvious structure of Hilbert $A$-$A$-bimodule, as any ideal of
$A$).}

\proofbegin Choosing a faithful representation as in 4.8 we may assume
that $A\subseteq \B(H_A)$, $B\subseteq \B(H_B)$ and $X\subseteq
\B(H_B, H_A)$. In addition the bimodule structure is composition of
operators and the inner products are given by $(x|y) = xy^*$ and $\<x,
y\> = x^*y$. The tensor product $X\*_BX^*$ is moreover identified with
$\overline{XX^*}$ (where the last occurrence of ``$*$'' should be
interpreted as the usual operator involution). This said, $X\*_BX^* =
(A|A)$.  \proofend

\beginsection 5. $K$-Theory and Hilbert Bimodules

The stage is now set for the presentation of the main section of this
work.  Given \cstar-algebras $A$ and $B$ as well as a left-full
Hilbert $A$-$B$-bimodule $X$, we want to define a group homomorphism
$X_*\:K_0(A)\arw K_0(B)$ which will be proven to be an isomorphism if
$X$ is right-full as well.

\statement{5.1.~Definition}{Let $X$ be a left-full Hilbert
$A$-$B$-bimodule.  If $\alpha$ is an element of the group $K_0(A)$,
which we identify with $F(A)$ under the isomorphism of 3.17, let $T$
be a Fredholm operator representing $\alpha$ in the sense that
$\ind(T) = \alpha$. We denote by $X_*(\alpha)$ the element of $K_0(B)$
defined by
  $$X_*(\alpha) = \ind(T\*I_X).$$}

Clearly, by 4.4, $X_*$ is well defined and it is not hard to see that
it is a group homomorphism.

\statement{5.2.~Proposition}{Let $X$ be a left-full Hilbert
$A$-$B$-bimodule and $Y$ be a left-full Hilbert $B$-$C$-bimodule, then
$Y_*\circ X_* = (X\*_B Y)_*$.}

\proofbegin This follows from the easy fact that tensor products are
associative, even if we drag along all the extra structure of Hilbert
bimodules. \proofend

Recall that the \cstar-algebras $A$ and $B$ are said to be strongly
Morita equivalent if there exists a Hilbert $A$-$B$-bimodule which is
simultaneously left and right-full. Such a module is called an
imprimitivity bimodule.  Given an imprimitivity bimodule $X$, we may
therefore consider the homomorphisms $X_*\:K_0(A) \arw K_0(B)$ and
$(X^*)_*\:K_0(B) \arw K_0(A)$, which compose to the identity in either
order by 5.2 and 4.13. The immediate outcome of these facts is our
main result.

\statement{5.3.~Theorem}{If $A$ and $B$ are strongly Morita equivalent
and $X$ is an imprimitivity bimodule, then $X_*$ is an isomorphism
{}from $K_0(A)$ onto $K_0(B)$.}

In order to treat $K_1$-groups we shall use the usual argument of
taking suspensions. The formalism of Hilbert bimodules developed in
section 4 makes is very easy to discuss Hilbert bimodules over tensor
product \cstar-algebras. So, before we embark on a study of $K_1$, let
us briefly deal with ``external tensor products''.

Let $A_i$ and $B_i$ be \cstar-algebras and $X_i$ be Hilbert
$A_i$-$B_i$-bimodules for $i=1, 2$.  Under faithful representations we
may assume that $A_i\subseteq \B(H_{A_i})$, $B\subseteq \B(H_{B_i})$
and $X_i\subseteq \B(H_{B_i}, H_{A_i})$.

Denote by $X_1\*X_2$ the closed linear span of the set
  $\{x_1\*x_2\:x_1\in X_1, x_2\in X_2\}$
  of operators in
  $\B( H_{B_1}\*H_{B_2} , H_{A_1}\*H_{A_2} )$.
   It is not hard to see that $X_1\*X_2$ is a Hilbert
$(A_1\*A_2)$-$(B_1\*B_2)$-bimodule. Here $A_1\*A_2$ and $B_1\*B_2$
mean the spatial tensor products.  If $x_i, y_i \in X_i$ for $i=1, 2$,
observe that
  $$(x_1\*y_1|x_2\*y_2) = (x_1\*y_1)(x_2\*y_2)^*
  = (x_1x_2^*)\*(y_1y_2^*) = (x_1|x_2)\*(y_1|y_2)$$
  so we see that the concrete inner-products on $X_1\*X_2$ may be
defined abstractly, at least on the algebraic tensor product of $X_1$
by $X_2$, without mentioning the representations. In fact it is easy
to see that the definition of $X_1\*X_2$ given above does not depend
(up to the obvious notion of isomorphism) on the particular faithful
representation chosen.

The object so defined is called the ``external'' tensor product of
$X_1$ and $X_2$ (compare \cite{\Jensen}, 1.2.4). It is readily
apparent that if each $X_i$ is left-full (resp.~right-full) then so is
$X_1\*X_2$.  As an obvious consequence we have:

\statement{5.4.~Proposition}{If $A_i$ and $B_i$ are \cstar-algebras
and if $A_i$ is strongly Morita equivalent to $B_i$ under the
imprimitivity bimodule $X_i$, for $i=1, 2$, then $A_1\*B_1$ is
strongly Morita equivalent to $A_2\*B_2$under the imprimitivity
bimodule $X_1\*X_2$.}

Let $C_0(\R)$ denote the \cstar-algebra of all continuous complex
valued functions on the real line. If $A$ and $B$ are strongly Morita
equivalent under the imprimitivity bimodule $X$, it follows that the
suspension of $A$, namely $SA = C_0(\R)\*A$ and and $SB$ are strongly
Morita equivalent to each other under the imprimitivity bimodule
$C_0(\R)\*X$. Using the standard isomorphism between $K_1(A)$ and
$K_0(SA)$, we have:

\statement{5.5.~Theorem}{If $A$ and $B$ are strongly Morita equivalent
and $X$ is an imprimitivity bimodule, then $(C_0(\R)\*X)_*$ is an
isomorphism from $K_1(A)$ onto $K_1(B)$.}

  \bigbreak
  \centerline{\tensc References}
  \nobreak\medskip
  \frenchspacing

  \bib{\Blackadar}
  {B. Blackadar}
  {$K$-theory for operator algebras}
  {MSRI Publications, Springer--Verlag, 1986.}

  \stdbib{\Brown}
  {L. G. Brown}
  {Stable isomorphism of hereditary subalgebras of \cstar-algebras}
  {Pacific J. Math.} {71} {1977} {335} {348}

  \stdbib{\BGR}
  {L. G. Brown, P. Green and M. A. Rieffel}
  {Stable isomorphism and strong Morita equivalence of
\cstar-algebras}
  {Pacific J. Math.} {71} {1977} {349} {363}

  \bib{\BMS}
  {L. G. Brown, J. A. Mingo and N. T. Shen}
  {Quasi-multipliers and embeddings of Hilbert \cstar-bimodules}
  {preprint, Queen's University, 1992}

  \stdbib{\Connes}
  {A. Connes}
  {Non-commutative differential geometry}
  {Publ. Math. IHES} {62} {1986} {257} {360}

  \bib{\Jensen}
  {K. Jensen and K. Thomsen}
  {Elements of \KK-Theory}
  {Birkh\"auser, 1991}

  \stdbib{\KI}
  {G. G. Kasparov}
  {Hilbert \cstar-modules: Theorems of Stinespring and Voiculescu}
  {J. Operator Theory} {4} {1980} {113} {150}

  \stdbib{\KII}
  {G. G. Kasparov}
  {The operator $K$-functor and extensions of \cstar-algebras}
  {Math. USSR Izvestija} {16} {1981} {513} {572}

  \stdbib{\Paschke}
  {W. Paschke}
  {Inner product modules over $B^*$-algebras}
  {Trans. Amer. Math. Soc.} {182} {1973} {443} {468}

  \stdbib{\RI}
  {M. A. Rieffel}
  {Induced representations of \cstar-algebras}
  {Adv. Math.} {13} {1974} {176} {257}

  \stdbib{\RII}
  {\sameauthor}
  {Morita equivalence for \cstar-algebras and $W^*$-algebras}
  {J. Pure Appl. Algebra} {5} {1974} {51} {96}

  \stdbib{\RIII}
  {\sameauthor}
  {Strong Morita equivalence of certain transformation group
\cstar-alge\-bras}
  {Math. Ann.} {222} {1976} {7} {22}

  \bib{\Rowen}
  {L. H. Rowen}
  {Ring theory -- Student edition}
  {Academic Press, 1991}

  \end